\documentclass{ifacconf}

\usepackage{graphicx}
\usepackage{natbib}        
\usepackage{subfigure}
\usepackage{amsmath}
\usepackage{bm}
\usepackage{amsfonts}
\usepackage{amssymb}
\usepackage{amsxtra}
\usepackage{cuted}
\usepackage{blindtext}
\usepackage{url}
\usepackage{suffix}
\usepackage{mathtools}
\usepackage{tabu}
\usepackage{caption}
\usepackage{adjustbox}
\usepackage{lipsum}
\usepackage{stfloats}
\usepackage{multirow}
\usepackage{tikz,colortbl}
\usepackage{zref-savepos}
\usepackage{pifont}
\usepackage{multicol}
\usepackage{tabu}
\usepackage{adjustbox}
\usepackage{multirow}
\usepackage{flushend}
\usepackage{float}
\usepackage{caption}
\usepackage{calc}
\usepackage{cancel}
\usepackage{mathrsfs}
\usepackage{algorithm}
\usepackage{algorithmic}
\usepackage{soul}
\usepackage{color}

\DeclareMathOperator*{\argmax}{arg\,max}
\DeclareMathOperator*{\argmin}{arg\,min}

\newtheorem{definition}{Definition}

\newtheorem{problem}{Problem}

\usetikzlibrary{calc}

\allowdisplaybreaks

\DeclarePairedDelimiter{\abs}{\lvert}{\rvert}

\newcommand{\norm}[1]{\left\lVert#1\right\rVert}

\begin{document}\sloppy
\begin{frontmatter}

\title{PolyAR: A Highly Parallelizable Solver For Polynomial Inequality Constraints Using Convex Abstraction Refinement} 

\thanks[footnoteinfo]{This work was partially sponsored by the NSF awards \#CNS-2002405 and \#CNS-2013824.}

\author[First]{Wael Fatnassi} 
\author[First]{Yasser Shoukry}

\address[First]{University of California, Irvine, California, USA (e-mail:\{wfatnass,yshoukry\}@uci.edu).}

\begin{abstract}                
: Numerical tools for constraints solving are a cornerstone to control verification problems. This is evident by the plethora of research that uses tools like linear and convex programming for the design of control systems. Nevertheless, the capability of linear and convex programming is limited and is not adequate to reason about general nonlinear polynomials constraints that arise naturally in the design of nonlinear systems. This limitation calls for new solvers that are capable of utilizing the power of linear and convex programming to reason about general multivariate polynomials. In this paper, we propose PolyAR, a highly parallelizable solver for polynomial inequality constraints. PolyAR provides several key contributions. First, it uses convex relaxations of the problem to accelerate the process of finding a solution to the set of the non-convex multivariate polynomials. Second, it utilizes an iterative convex abstraction refinement process which aims to prune the search space and identify regions for which the convex relaxation fails to solve the problem. Third, it allows for a highly parallelizable usage of off-the-shelf solvers to analyze the regions in which the convex relaxation failed to provide solutions. We compared the scalability of PolyAR against Z3 8.9 and Yices 2.6 on control designing problems. Finally, we demonstrate the performance of PolyAR on designing switching signals for continuous-time linear switching systems.
\end{abstract}


\end{frontmatter}
\section{INTRODUCTION}
\vspace{-3mm}
Advances in constraints programming have opened several venues for control system synthesis and verification of hybrid systems. For instance, linear programming and convex optimization are heavily used in a multitude of control system design and analysis tools. Recent surveys [\cite{surveuNumTool1,surveuNumTool2}] showed that such numerical tools had changed the control system design philosophy.

Nevertheless, linear and convex programming are limited in their ability to problems with specific structures. In several hybrid system design and verification problems, constraints are neither linear nor convex. This calls for efficient solvers that can reason about \textit{general multivariate polynomial constraints}. In that regard, Cylindrical Algebraic Decomposition (CAD) has long been one of the most influential algorithms capable of solving general multivariate polynomial constraints. The first CAD algorithm was introduced by [\cite{collins}]. However, the use of CAD is often limited by the number of variables in the input polynomials, a reflection of its worst-case complexity that grows in a doubly exponential fashion in the number of variables [\cite{complexityproblem1}].

To alleviate the CAD's doubly exponential issue, we introduce PolyAR, a highly parallelizable solver that uses convex programming and abstraction refinement to solve general multivariate polynomial inequality constraints. The main novel contributions of this work can be summarized as follows:
\begin{itemize}
    \item PolyAR is a highly parallelizable solver that uses a combination of convex programming and abstraction refinement to solve multivariate polynomial inequality constraints.  
    \item PolyAR uses a novel convex abstraction refinement process where the original problem is iteratively relaxed into a series of convex programming problems with the aim to find the solution and prune the search space. Second, it refines such abstraction where it becomes tighter with each iteration of the algorithm. Finally, it examines in parallel all the identified small volume regions left from the abstraction using off-the-shelf solvers (e.g., Z3 and Yices) to search for a solution in these regions. 
    \item We validate our approach by comparing the scalability of the proposed PolyAR solver with respect to the latest versions of state-of-art non-linear real arithmetic solvers, such as Z3 8.9 and Yices 2.6, on synthesizing stabilizing static output feedback controller (SOF) for linear time-invariant (LTI) continuous systems and designing non-parametric controller for the non-linear Duffing oscillator.
    \item We demonstrate the performance of PolyAR on the problem of designing switching signals for continuous-time linear switching systems.
\end{itemize}
\vspace{-1mm}

\textbf{Related work:} The original CAD algorithm that was introduced by \cite{collins} was the first algorithm that solves general polynomial inequality constraints. However, due to Collins CAD's high time complexity, there have been improvements to this algorithm. \cite{Hong} proposed an improvement of the projection operator in Collins CAD. However, the execution time of the modified version of the algorithm is still limited by the number of variables. \cite{McCalum} introduced a new projection operator which is a subset of Hong projection operator, removing redundant polynomials. However, McCallum proved that lifting over a sign-invariant CAD with this projection set is not sufficient to guarantee sign-invariance which makes the algorithm prone to error. The ABsolver tool proposed by \cite{ABsolver} leverages a generic nonlinear optimization tool for solving non-linear constraints. However, generic optimization tool may produce incomplete results, and possibly incorrect, due to the local nature of the solver. \cite{Z3} introduced Z3 which is another solver that implements an efficient nonlinear real arithmetic solver, that provide support for nonlinear polynomial arithmetic. However, it is still affected a lot by the increase in the number of variables in the polynomials. Because of the high complexity of existing approaches, we propose a highly parallelizable, efficient, and complete solver that uses the advantage and the simplicity of convex optimizations and abstraction refinement to solve higher order polynomial inequality constraints. To the best of our knowledge, this approach is new and has not been highlighted before.


\section{Problem Formulation}
\vspace{-3mm}
\subsection{Notation}
\vspace{-3mm}
We denote by $x=\big(x_1,x_2,\cdots,x_n\big) \in \mathbb{R}^n$ the set of real-valued variables, where $x_i \in \mathbb{R}$. We denote by $I_n=\big[\underline{d}_1,\overline{d}_1\big] \times \cdots \times$ $\big[\underline{d}_n,\overline{d}_n\big] \subset \mathbb{R}^{n}$ the $n$-dimensional region. We denote the space of polynomials with $n$ variables and coefficients in $\mathbb{R}$ by $\mathbb{R}[\left(x_1,x_2,\cdots,x_n\right)]$.  We denote by $\wedge$ the Boolean conjunction. A set of the form $L^{-}_{0}\left(f\right)=\{\big(x_1,\cdots,x_n\big)\big|f\big(x_1,\cdots,x_n\big)\leq 0\}$ ($L^{+}_{0}\left(f\right)=\{\big(x_1,\cdots,x_n\big)\big|f\big(x_1,\cdots,x_n\big)\geq 0\}$) is called zero sublevel (superlevel) set of $f$, respectively.

\subsection{Main Problem}
In this paper, we focus on \textit{polynomial inequality constraints} with input ranges as closed boxes which are described in the following definition:
\begin{definition}
A polynomial inequality constraint $F=I^n~\wedge~P_m$ consists of:
\begin{itemize}
    \item a set of interval constraints:
    \begin{align}\label{pic1}
      I^n=\bigwedge\limits_{i=1}^{n}x_i\in[\underline{d}_i,\overline{d}_i],  
    \end{align}
    \item a polynomial constraint:
    \begin{align}\label{pic2}
     P_m=\bigwedge\limits_{i=1}^{m}~p_i\left(x_1,\cdots,x_n\right)~\leq~0,   
    \end{align}
\end{itemize}
where $p_i\left(x\right)=p_i\left(x_1,\cdots,x_n\right) \in \mathbb{R}[\left(x_1,x_2,\cdots,x_n\right)]$ is a polynomial over variables  $x_1,\cdots,x_n$. Without loss of generality, $\bigwedge\limits_{i=1}^{m}~p_i\big(x\big)$ $\geq~0$ and $\bigwedge\limits_{i=1}^{m}~p_i\left(x\right)~=~0$ can be encoded in constraint number \eqref{pic2}.
\end{definition}

We are now in a position to state the problem that we will consider in this paper.

\begin{problem}
$\exists x=\left(x_1,\cdots,x_n\right)$ subject to $F=I^n~\wedge~P_m$. 
\end{problem}


\section{Abstraction Refinement of Higher Order Polynomials Using Quadratic Polynomials} \label{sec:abstraction}
\vspace{-3mm}
Traditional techniques for solving Problem 1 focus on finding all the $n$ roots of the $m$ polynomials and check all the regions between two successive roots to assign a positive/negative sign for each of these regions. Therefore, solving Problem 1 is known to be a doubly combinatorial problem in $n$ with a total running time that is bounded by $\left(md\right)^{2^n}$ \cite{complexityproblem1}, where $d$ is the maximum degree among polynomials in $P_m$.

In problems that are doubly exponential in the input space $n$, it is beneficial to isolate subsets of the search space in which the solution is guaranteed not to exist. Recall that Problem 1 asks for an $x$ in $\mathbb{R}^n$ for which all the polynomials are negative. Therefore, a solution does not exist in subsets of $\mathbb{R}^n$ at which one of the polynomials is always positive. Similarly, isolating regions of the input space for which some of the polynomials are negative is also beneficial to finding the solution faster.

Our tool's main novelty is to use ``convex abstractions'' of the polynomials to find subsets of $L^{+}_{0}\left(p\right)$ and $L^{-}_{0}\left(p\right)$ efficiently. Indeed such ``abstractions'' may not be able to identify all regions for which the polynomial is positive or negative, which calls for an ``abstraction refinement'' process in which these ``convex abstractions'' become tighter with each iteration of the algorithm.

\begin{figure}
	\includegraphics[width=0.9\columnwidth, trim=0 50mm 0 0, clip]{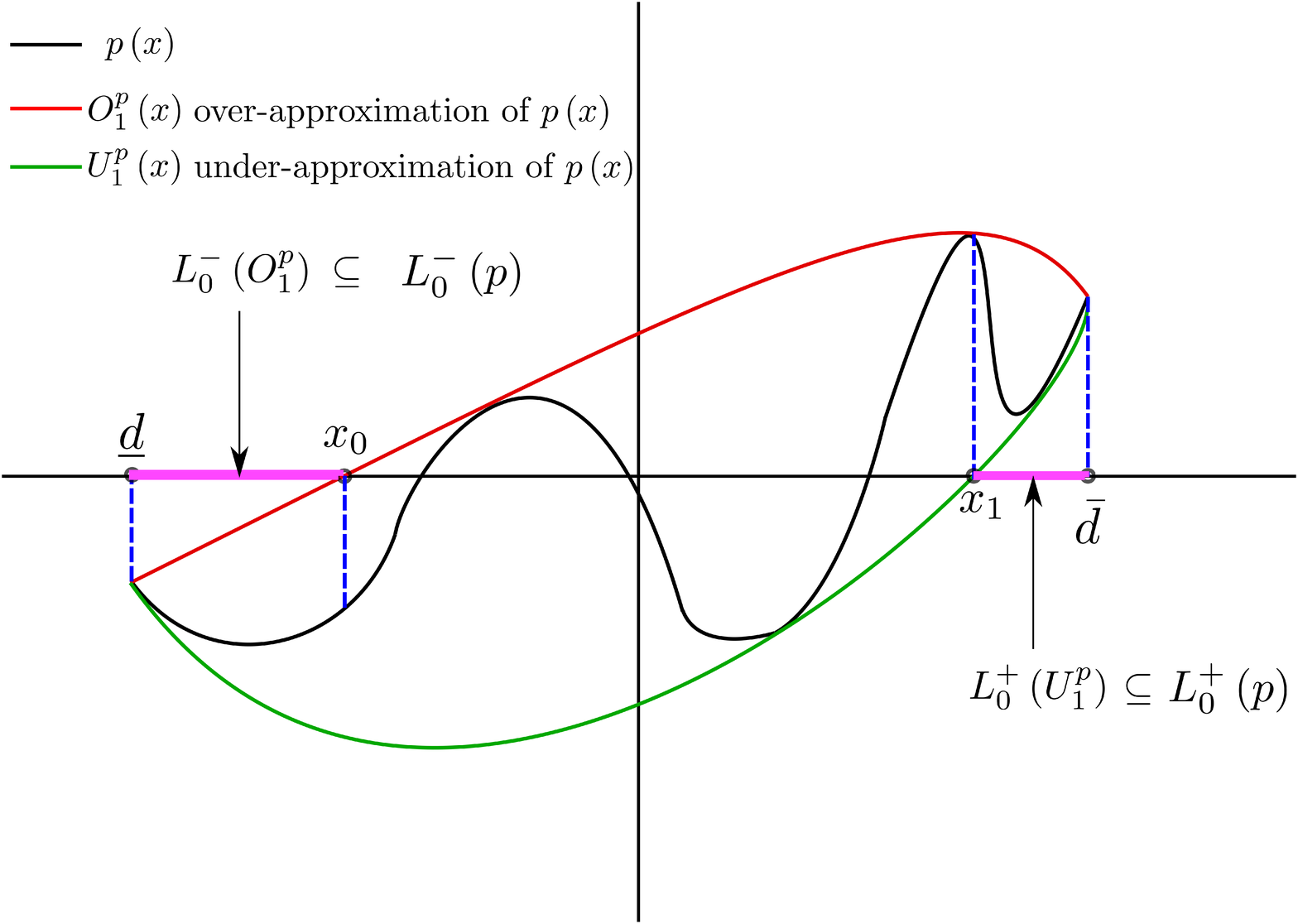} \\
	\includegraphics[width=0.9\columnwidth, trim=0 70mm 0 0, clip]{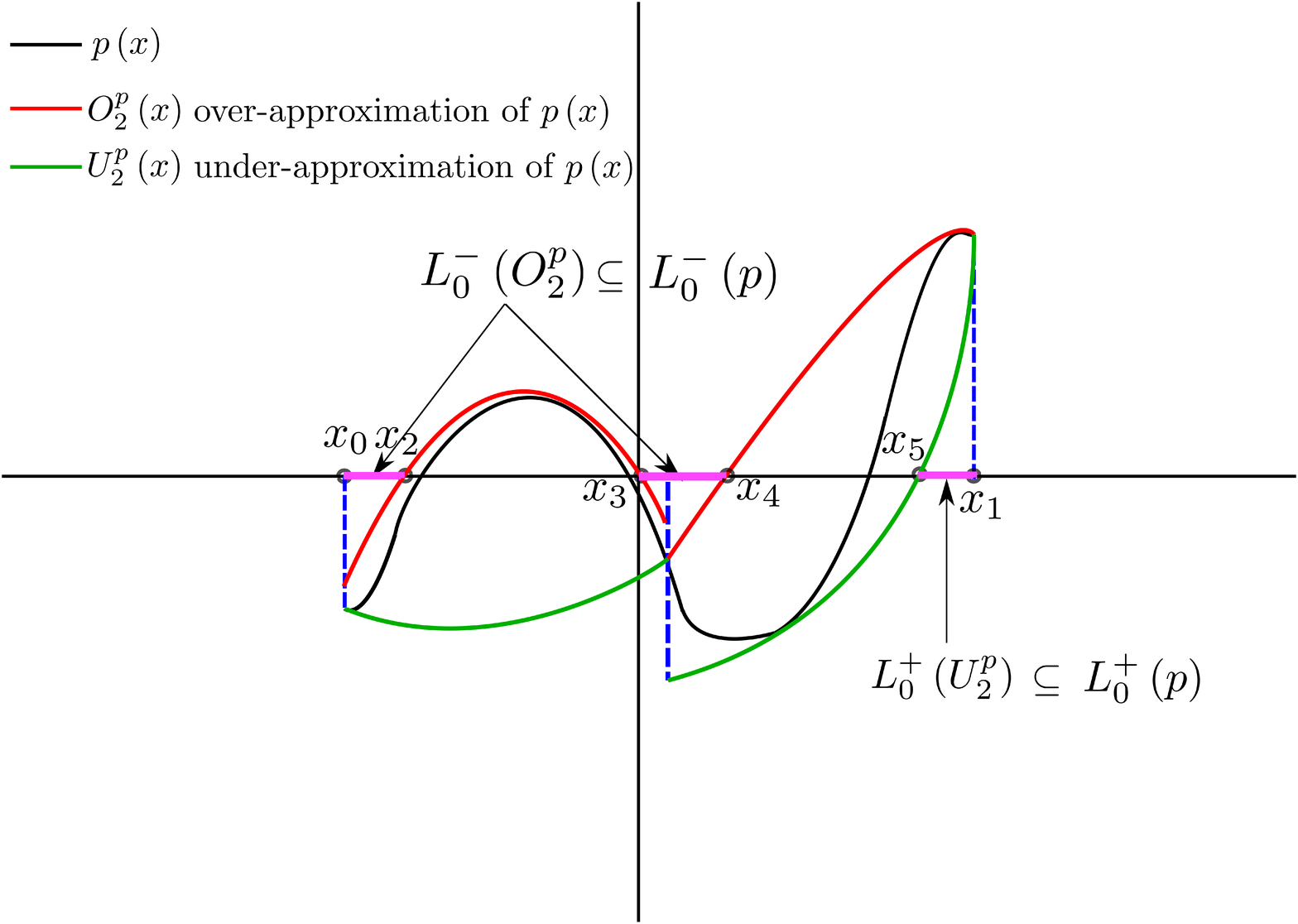}
	\caption{Abstraction Refinement of higher order polynomial using quadratic approximations: (top) first iteration and (bottom) second iteration.}
	\label{F1}
\end{figure}

Figure~\ref{F1}(top) visualizes the proposed abstraction refinement process. Starting from a polynomial $p(x) \in \mathbb{R}[x]$ and an interval $I_n \subset \mathbb{R}^n$, we compute two quadratic polynomials:
\begin{align*}
    O^p_1(x) &\ge p(x) \qquad \forall x \in I_n,\\
    U^p_1(x) &\le p(x) \qquad \forall x \in I_n.
\end{align*}
where $O$ and $U$ stands for Over-approximate and Under-approximate quadratic polynomials, respectively, and the subscript in $O^p_1(x)$ and $U^p_1(x)$ encodes the iteration index of the abstraction refinement process. Computing such upper and lower abstractions can be carried out efficiently using Taylor approximation. Please refer to the example depicted in Figure~\ref{F1} (top) for a visualization of $O^p_1(x)$ and $U^p_1(x)$ for one dimensional higher order polynomial (order $\geq$ 3)  defined in the closed interval $[\underline{d},~\overline{d}]~\subset~\mathbb{R}$.

The next step is to use the quadratic abstractions to isolate subsets of $L^{+}_{0}\left(p\right)$ and $L^{-}_{0}\left(p\right)$. It is particularly direct to show that the zero superlevel set of $U^p_1(x)$ is a subset of $L^{+}_{0}\left(p\right)$, i.e., $L^{+}_{0}\left(U^p_1\right) \subseteq L^{+}_{0}\left(p\right)$. Similarly, the zero sublevel set of $O^p_1(x)$ is a subset of $L^{-}_{0}\left(p\right)$, i.e., $L^{-}_{0}\left(O^p_1\right) \subseteq L^{-}_{0}\left(p\right)$. Thanks to the fact that $O^p_1(x)$ and $U^p_1(x)$ are quadratic polynomials, finding their zero superlevel and zero sublevel sets, respectively, can be computed efficiently. Referring to the example in Figure~\ref{F1}(top), these zero superlevel and sublevel sets are $L^{+}_{0}\left(U^p_1\right)=[x_1, \overline{d}]$ and $L^{-}_{0}\left(O^p_1\right)=[\underline{d},x_0]$, respectively. 

It is clear from Figure~\ref{F1}(top) that the abstractions $O^p_1(x)$ and $U^p_1(x)$ fails to identify all subsets of $L^{-}_{0}\left(p\right)$ and $L^{+}_{0}\left(p\right)$. Therefore, the next step is to compute tighter over and under approximations of $p(x)$. Such a refinement process can be carried out by removing the zero superlevel and the zero sublevel sets, i.e., $L^{+}_{0}\left(U^p_1\right)$ and $L^{-}_{0}\left(O^p_1\right)$, identified using the previous abstraction and computing new over and lower approximation, as shown in Figure~\ref{F1}(bottom). The process of abstraction refinement can continue until the remaining subsets of the search space, in which case we call them ambiguous regions, and with some abuse of notation, denoted them by $L^{+/-}_{0}\left(p\right)$, are \emph{small} enough to be analyzed using off-the-shelf solvers. More details about the proposed abstraction refinement process are given in the next section.

\vspace{-3mm}
\section{Algorithm Architecture}
\vspace{-3mm}
\begin{figure}[!t]
	\includegraphics[width=1.0\columnwidth]{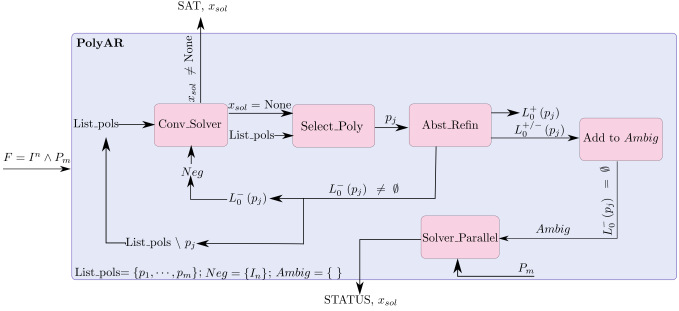} 
	\caption{Framework of PolyAR.} 
	\label{fig:PolyAR}
\end{figure}
In this section, we describe the different steps used by our solver PolyAR to solve Problem 1. 

Our design methodology for the PolyAR tool aims to reduce the number of the required abstraction refinement and tries to find a solution early on in the process. To that end, the tool starts by computing a set of convex (quadratic or linear) polynomials $O_0^{p_i}(x), i = 1,\ldots,m$, that over approximate the original polynomials. The next step is to solve a convex feasibility problem aiming to find a solution that satisfy the constraints:
$$ \exists x = (x_1,\ldots,x_n) \qquad \text{s.t.} \qquad  O_0^{p_i}(x) \le 0, \quad i = 1,\ldots,m.$$
Indeed, if such a convex problem is feasible, the tool terminates and returns the solution found by the convex feasibility problem above (\textbf{Conv\_Solver}, Line~\ref{alg:line:conv} in Algorithm~\ref{alg:PolyAR}). If not, then the tool selects one polynomial $p_j$ (\textbf{Select\_Poly}, Line \ref{alg:line:selectPoly}) to perform the abstraction refinement process. Indeed, several heuristics can be applied to select which polynomial will be selected. In the PolyAR tool, we opt-out to select the polynomial with the highest Lipschitz constant. Our intuition is that the higher the Lipschitz constant, the harder to obtain a tight over-approximation that can be used to find the solution.

Once a polynomial $p_j$ is selected, the next step is to apply the abstraction refinement process on $p_j$  (\textbf{Abst\_Refin}, Line~\ref{alg:line:abst_ref} in Algorithm~\ref{alg:PolyAR}). The objective of the abstraction refinement process is to identify subsets of the positive regions $L^{+}_{0}\left(p_j\right)$  and negative regions $L^{-}_{0}\left(p_j\right)$. Indeed, such abstraction refinement may not be able to identify all positive and negative regions, and hence a remaining portion of the search space may not be identified to belong to either $L^{+}_{0}\left(p_j\right)$ or $L^{-}_{0}\left(p_j\right)$ in which case it belongs to the ambiguous region $L^{+/-}_{0}\left(p_j\right)$. The abstraction refinement process of the polynomial $p_j$ ensure that the volume of such ambiguous regions are below a certain user defined threshold.

The process of using the convex solver to find the solution and abstracting one polynomial continues. Since a solution of Problem 1 needs to lie in a negative region for all the polynomials, we confine the tool attention to the negative regions identified by the abstraction refinement in the previous iterations (Line~\ref{alg:line:neg} in Algorithm~\ref{alg:PolyAR}) to accelerate the process of searching for the solution. 

While excluding the positive regions identified in previous iterations does not affect the tool (since a solution is guaranteed not to exist in such regions), excluding the ambiguous regions from the next iterations may affect the correctness of the tool. Therefore, the last step in the PolyAR tool is to examine all the identified ambiguous regions using off-the-shelf solvers (e.g., Z3 and Yices) to search for a solution in these regions (\textbf{Solver\_Parallel}, Line~\ref{alg:line:solver} in Algorithm~\ref{alg:PolyAR}). Because the volume of these ambiguous regions is smaller than a user-defined threshold, the execution time of running off-the-shelf tools on such small volume regions is shorter than solving the original problem. This reflects that the number of roots for each polynomial is limited in small regions. Moreover, searching for a solution in these ambiguous regions can be highly parallelized, leading to an extra level of efficiency. This process is summarized in Algorithm 1 and Figure~\ref{fig:PolyAR}. We describe in detail each block algorithm that constitutes Algorithm 1 in the next subsections.

\begin{algorithm}
\caption{PolyAR$\left(F\right)$}  \label{alg:PolyAR}
\begin{flushleft}
\textbf{Input:} $F = I^n \wedge P_m$\\
\textbf{Output}: STATUS, $x_{\text{Sol}}$
\end{flushleft}
\begin{algorithmic}[1]
\STATE $Neg= \{I_n \}$ 
\STATE $Ambig= \{\}$
\STATE $\text{List\_pols}= \{p_1,\ldots,p_m\}$
\WHILE{$\text{List\_pols}~\neq~ \emptyset$}
    \STATE $x_{\text{Sol}}:=\textbf{Conv\_Solver}\left(Neg,\text{List\_pols}\right)$ \label{alg:line:conv}
    \IF{$x_{\text{Sol}}~\neq~\text{None}$}
    \STATE STATUS=SAT
    \RETURN STATUS, $x_{\text{Sol}}$
    \ENDIF

     \STATE $p_j=\textbf{Select\_Poly}\left(\text{List\_pols}\right)$ \label{alg:line:selectPoly}
     \STATE $L^{-}_{0}\left(p_j\right),L^{+}_{0}\left(p_j\right),L^{+/-}_{0}\left(p_j\right)$ \\ \quad \quad \quad \quad \quad \quad \quad \quad \quad \quad \quad $:= \textbf{Abst\_Refin}\left(Neg,p_j\right)$ \label{alg:line:abst_ref}
      \STATE $Ambig.\text{add}\left(L^{+/-}_{0}\left(p_j\right)\right)$    
     \IF{$L^{-}_{0}\left(p_j\right)~==~\emptyset$}
         \STATE break
     \ENDIF
     
     \STATE $Neg=L^{-}_{0}\left(p_j\right)$ \label{alg:line:neg}
     \STATE $\text{List\_pols}=\text{List\_pols}\setminus{p_j}$
\ENDWHILE
\IF{$\text{List\_pols}~\neq~ \emptyset$}
    \STATE $\text{STATUS}, x_{\text{Sol}}:=\textbf{Solver\_Parallel}\left(Ambig,P_m\right)$ \label{alg:line:solver}
         
         \RETURN $\text{STATUS}, x_{\text{Sol}}$
\ELSE
    \STATE $\text{STATUS=SAT}$
    \STATE $x_{\text{Sol}}=\text{center}\left(Neg\right)$
    \RETURN $\text{STATUS}, x_{\text{Sol}}$
\ENDIF    
\end{algorithmic}
\end{algorithm}

\vspace{-2mm}
\subsection{Early Termination Using \textbf{Conv\_Solver}:}
The objective of the \textbf{Conv\_Solver} (Algorithm 2) is to search for a solution to Problem 1 using the information of (i) a set of closed convex regions $Neg$ identified by the previous iterations of the abstraction refinement process and (ii) a list of polynomials (List\_pols) that have not yet been processed by the abstraction refinement process.

Our approach is to compute a convex over-approximation of the polynomials in List\_pols using Taylor approximation. To that end, we recall the definition of Taylor polynomials: 
\begin{definition}
Let $f:\mathbb{R}^n\rightarrow\mathbb{R}$ be two times differentiable in open interval around a point $a \in \mathbb{R}^n$, then $f\left(x\right)$ can be written in terms of first and second order Taylor polynomials , $T_1\left(x\right)$ and $T_2\left(x\right)$, around the neighborhood of $a$, as follows:
\vspace{-2mm}
\begin{eqnarray}\label{eq6}
f\left(x\right) &=& \left(x-a\right)^{T}D_{f}\left(a\right)+R_1\left(c_1\right)\nonumber\\
&=&T_1\left(x\right)+R_1\left(c_1\right), \\
f\left(x\right) &=& \left(x-a\right)^{T}D_{f}\left(a\right)\nonumber\\
&+&\frac{1}{2}\left(x-a\right)^{T}H_{f}\left(a\right)\left(x-a\right)+R_2\left(c_2\right)\nonumber\\
&=&T_2\left(x\right)+R_2\left(c_2\right),
\end{eqnarray}
\end{definition}
\vspace{-2mm}
\noindent where $T_1\left(x\right)=\left(x-a\right)^{T}D_{f}\left(a\right)$ and $T_2\left(x\right)=\left(x-a\right)^{T}D_{f}\left(a\right)+\frac{1}{2}\left(x-a\right)^{T}H_{f}\left(a\right)\left(x-a\right)$. $D_{f}\left(a\right)$ and $H_{f}\left(a\right)$ denote the Gradient vector and the Hessian matrix of $f$ at the point $a$. $R_1\left(c_1\right)$ and $R_2\left(c_2\right)$ are reminders that depend on $a$ and two points $c_1 \in \mathbb{R}^n$ and $c_2 \in \mathbb{R}^n$ that are located in the neighborhood of $a$.
We upper-bound $R_1\left(c_1\right)$ and $R_2\left(c_2\right)$ by $M_1 \in \mathbb{R}_{+}$ and $M_2 \in \mathbb{R}_{+}$, i.e., $\abs{R_1\left(c_1\right)} \leq M_1$ and $\abs{R_2\left(c_2\right)} \leq M_2$. Hence:
\begin{align*}
    T_1(x) - M_1 &\le f(x) \le T_1(x) + M_1 \\
    T_2(x) - M_2 &\le f(x) \le T_2(x) + M_2
\end{align*}
Next, we check the convexity of the obtained second Taylor approximation and use it to compute the over-approximation function $O^{p_i}$ whenever it is convex (Line~\ref{line:alg:taylor2} in Algorithm~\ref{alg:conv}). Otherwise, we use the first Taylor approximation instead (Line~\ref{line:alg:taylor1} in Algorithm~\ref{alg:conv}). Finally, for each $region$ in the set of negative regions ($Neg$), we solve the following convex feasibility problem:
\begin{align}\label{eq10}
x_{\text{Sol}}:=&\argmin\limits_{x \in region}^{}1 ~~ \text{s.t.} ~~ O^{p_i}\left(x\right) \leq~0,~~i \in \text{List\_pols}.&
\end{align}

\begin{algorithm}[t]
\caption{$\textbf{Conv\_Solver}\left(Neg,\text{List\_pols}\right)$} \label{alg:conv}
\begin{flushleft}
\textbf{Input:} $Neg$, $\text{List\_pols}$\\
\textbf{Output}: $x_{\text{Sol}}$
\end{flushleft}
\begin{algorithmic}[1]
\FOR{$region \in Neg$}
\FOR{$i \in \text{List\_pols}$}
    \IF{$\texttt{Taylor}_{over}\left(p_i\left(x\right), 2\right)~\text{is convex}$}
         \STATE $O^{p_i}\left(x\right)=\texttt{Taylor}_{over}\left(p_i\left(x\right), 2\right)$ \label{line:alg:taylor2}
        
    \ELSE
         \STATE $O^{p_i}\left(x\right)=\texttt{Taylor}_{over}\left(p_i\left(x\right), 1\right)$ \label{line:alg:taylor1}
    \ENDIF
\ENDFOR

\STATE $x_{\text{Sol}}:=\argmin \limits_{x \in region}^{}1 \quad s.t. \quad O^{p_i}\left(x\right)\leq~0$ (see eq. \eqref{eq10}))
\RETURN $x_{\text{Sol}}$
\ENDFOR
\end{algorithmic}
\end{algorithm}
\vspace{-5mm}
\subsection{Abstraction Refinement Using \textbf{Abst\_Refin}:}
Given the set of negative regions $Neg$ identified by the previous abstraction refinement process along with the polynomial $p_j$ selected by the \textbf{Select\_Poly} algorithm, the objective of the \textbf{Abst\_Refin} algorithm is to find subsets of the zero sublevel sets of $p_j$ that lie inside $Neg$. The output of this algorithm are subsets of $L_0^-(p_j)$ and $L_0^+(p_j)$. The remainder of $Neg$ is then considered to be part of the ambiguous regions $L_0^{+/-}(p_j)$. To do so, for every $region$ in $Neg$, the tool initiates a list of ambiguous regions $List\_Ambig\_reg$, which will contain all the ambiguous regions from the abstraction refinement (Line 4 in Algorithm 3). Next, it selects one element from these ambiguous regions (Line 5 in Algorithm 3) and performs the abstraction refinement on this region iteratively until the volume of the remaining ambiguous region is smaller than a user-defined threshold (Line 6 in Algorithm 3). During the iterative abstraction refinement, all the identified zero sublevel and superlevel subsets are stored in the sets $L_0^-(p_j)$ and $L_0^+(p_j)$, respectively.

\begin{figure}[!ht]
		\centering
		\includegraphics[width=0.7\columnwidth]{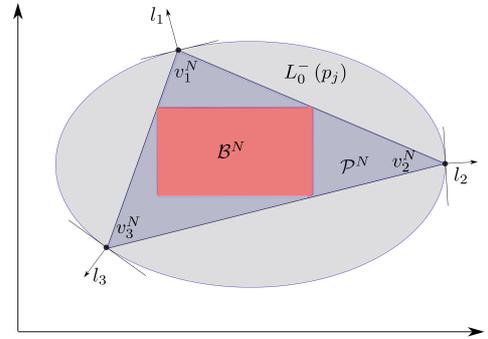}
		\caption{\small{Polytopic under-approximation of a $2-$dimensional ellipse sublevel set $L^{-}_{0}\left(p_j\right)$. $\mathcal{P}^{N}$ presents the under-approximate polytope inscribed in $L^{-}_{0}\left(p_j\right)$, and $\mathcal{B}^{N}$ represents the axis-aligned box of maximum volume inscribed in $\mathcal{P}^{N}$.}}
		\label{undpol}
\end{figure}

While the zero sublevel (superlevel) sets of the quadratic over-approximation (under-approximation) are ellipsoid or hyperboloid in general, we opt to represent all the subsets of $L_0^-(p_j)$ and $L_0^+(p_j)$ as $n-$dimensional hypercubes. This choice reflects the fact that off-the-shelf solvers (e.g., Z3 and Yices) can exploit the geometry of hypercubes to accelerate their computations. The process of finding these hypercubes can be summarized as follows: 
\begin{enumerate}
    \item \textbf{Step 1:} Compute the largest polytope inside the ellipsoid or hyperboloid representing the zero sublevel (superlevel) sets of the quadratic over-approximation (under-approximation) of $p_j$. To that end, we use a set of user-defined templates for the polytope.
    \item \textbf{Step 2:} The previous step uses user-defined templates to find the polytope, such templates may fail and return an infeasible solution. In such scenarios, we split the ambiguous region into two (along the longest dimension) until a polytope is found.
    \item \textbf{Step 3:} Finally, we under approximate the computed polytope with hypercubes.
\end{enumerate}
This process is visualized in Figure~\ref{undpol}. The details of each of these steps are given in the following subsections.

\subsubsection{\textbf{Step 1: Computing the largest polytope subset of $L_0^-(p_j)$ and $L_0^+(p_j)$}}

Given the over-approximation $O^{p_j}$ computed using Taylor polynomials (detailed in Section 4.1) and a convex ambiguous region (Ambig\_reg), we start by computing a set of $n+1$ vertices $v^{N}_1, \ldots, v^{N}_{n+1}$ that are inscribed in the ambiguous region $Ambig\_reg$. Each vertex can be computed by solving the following convex optimization problem:
\begin{align}\label{vertn}
v^{N}_i=&\argmin \limits_{v_i \in Ambig\_reg} \left(l^T_iv_i\right) 
    \qquad  \text{s.t.} \qquad O^{p_j}(v_i)~\leq~0,& 
\end{align} 
where $l_i$ is a user defined normal vector (or template) (see Figure~\ref{undpol} for graphical representation of such normal vectors). Using these vertices, we can obtain the polytope $\mathcal{P}^{N}$ as:
$$\mathcal{P}^{N} = \textbf{Convex\_Hull}\left(v^{N}_1, \ldots,  v^{N}_{n+1} \right).$$
Thanks to the constraints in the optimization problem~\eqref{vertn} along with the convexity of $L_0^-(p_j)$, it is direct to conclude that the polytope $\mathcal{P}^{N}$ satisfy $\mathcal{P}^{N} \subset L_0^-(p_j)$. We compute the polytope $\mathcal{P}^{P} \subset L_0^+(p_j)$ in a similar fashion using the under-approximation $U^{p_j}$ (Line 23 in Algorithm 3).

\begin{algorithm}[t]
\caption{$\textbf{Abst\_Refin}\left(Neg,p_j\right)$}
\begin{flushleft}
\textbf{Input:} $Neg$, $p_j$ \\
\textbf{Output}: $L^{-}_{0}\left(p_j\right)$, $L^{+}_{0}\left(p_j\right)$, $L^{+/-}_{0}\left(p_j\right)$
\end{flushleft}
\begin{algorithmic}[1]
\STATE  $L^{-}_{0}\left(p_j\right)=\{~\}$, $L^{+}_{0}\left(p_j\right)=\{~\}$, $L^{+/-}_{0}\left(p_j\right)=\{~\}$
\FOR{$region \in Neg$}
\STATE $\text{vertices}^{N}=\{~\}$, $\text{vertices}^{P}=\{~\}$
\STATE $List\_Ambig\_reg=\{region\}$
\STATE $Ambig\_reg=\textbf{Select\_region}\left(List\_Ambig\_reg\right)$
\WHILE{$\text{Volume}\left(Ambig\_reg\right)~>~\text{Vol}_{\text{threshold}}$}
    \STATE $List\_Ambig\_reg=List\_Ambig\_reg\setminus{Ambig\_reg}$
    \FOR{$i~\in~\big(1,\cdots,n+1\big)$}
        \STATE  $v^{N}_i=\argmin\limits_{v_i \in Ambig\_reg} \left(l^T_iv_i\right) \quad s.t. \quad O^{p_j}(v_i) ~\leq~0.$
        \IF{$v^{N}_i~\neq~\text{None}$} 
              \STATE $\text{vertices}^{N}.\text{add}\left(v^{N}_i\right)$ 
        \ENDIF
        
        \STATE $v^{P}_i=\argmin\limits_{v_i \in Ambig\_reg} \left(l^T_iv_i\right) \quad s.t. \quad U^{p_j}(v_i)~\leq~0.$
        \IF{$v^{P}_i~\neq~\text{None}$} 
             \STATE $\text{vertices}^{P}.\text{add}\left(v^{P}_i\right)$  
        \ENDIF
    \ENDFOR
    
\IF{$\big(\text{vertices}^{N}==\emptyset ~ \text{and} ~ \text{vertices}^{P}==\emptyset\big)$}
    \STATE $Ambig\_reg_{1}$, $Ambig\_reg_{2}$ \\ \quad \quad \quad \quad \quad \quad \quad \quad \quad  $:=\text{\textbf{Half\_Div}}\left(Ambig\_reg\right)$
    \STATE $List\_Ambig\_reg.add\left(Ambig\_reg_{1}, Ambig\_reg_{2}\right)$

\ELSIF{$\big(\text{vertices}^{N}\neq \emptyset  ~ \text{and} ~ \text{vertices}^{P}\neq \emptyset\big)$}
    \STATE $\mathcal{P}^N=\textbf{Convex\_Hull}\left(\text{vertices}^{N}\right)$
    \STATE $\mathcal{P}^P=\textbf{Convex\_Hull}\left(\text{vertices}^{P}\right)$
     
    \STATE $\mathcal{B}^{N}=\textbf{Box}\left(\mathcal{P}^N\right)$; $\mathcal{B}^{P}=\textbf{Box}\left(\mathcal{P}^P\right)$ 
    \STATE $L^{-}_{0}\left(p_j\right).\text{add}\left(\mathcal{B}^{N}\right)$; $L^{+}_{0}\left(p_j\right).\text{add}\left(\mathcal{B}^{P}\right)$ 

    \STATE $Ambig\_reg=Ambig\_reg\setminus{\left(\mathcal{B}^{N} \cup \mathcal{B}^{N}\right)}$
    
    \STATE $List\_Ambig\_reg.\text{add}\big(Ambig\_reg\big)$

\ENDIF
\STATE $Ambig\_reg=\textbf{Select\_region}\left(List\_Ambig\_reg\right)$
\ENDWHILE
\STATE $L^{+/-}_{0}\left(p_j\right).add\left(List\_Ambig\_reg\right)$
\ENDFOR
\RETURN $L^{-}_{0}\left(p_j\right)$, $L^{+}_{0}\left(p_j\right)$, $L^{+/-}_{0}\left(p_j\right)$
\end{algorithmic}
\end{algorithm}

\subsubsection{\textbf{Step 3: Under approximate the polytopes with axis aligned boxes:}}
To compute the largest axis-aligned hypercube $\mathcal{B}^{N}$ inscribed inside the polytope $\mathcal{P}^{N}$, we solve the following convex optimization problem~\cite{InscribRec}:
\begin{align}\label{eq14}
& \argmax\limits_{(l_1^N, u_1^N, \ldots, l_n^N,u_n^N)\in \mathbb{R}^{2n}}^{} \sum\limits_{k=1}^{n}\log\left(u^N_k-l^N_k\right)\nonumber\\
& \text{s.t.}  \sum\limits_{k=1}^{n}\left(p^{N,+}_{ik}u^N_{k}-p^{N,-}_{ik}l^N_{k}\right)\leq c^N_i,~i=1,\cdots,n_p,
\end{align}
where $(l_1^N, u_1^N, \ldots, l_n^N,u_n^N) \in \mathbb{R}^{2n}$ is the representation of the box $\mathcal{B}^{N}$ with $(l_k,u_k)$ is the lower/upper limit of the box in the $k$th dimension, $p_{ik}^{N,+}=\max\{p^N_{ik},0\}$, $p_{ik}^{N,-}=\max\{-p^N_{ik},0\}$, and $p^N_{ik},c^N_i$ are the rows of the half-space matrix/vector representation of the polytope $\mathcal{P}^{N}$.

\subsection{Highly Parallelizable Analysis of Ambiguous Regions using \textbf{Solver\_Parallel}}
\vspace{-2mm}
Once all the ambiguous regions are identified, the next step is to analyze all of them using off-the-shelf solvers. In particular, PolyAR supports the use of the latest versions Z3 8.9 and Yices 2.6 solvers. Thanks to the fact that all the ambiguous regions are hypercubes, both these solvers can exploit the geometry of the region to accelerate their computations. Also, thanks to the fact that the volume of all ambiguous regions is lower than a user-defined threshold, the CAD algorithm can run efficiently. To that end, PolyAR tool runs multiple instances of Z3 or Yices to analyze all these ambiguous regions in parallel as summarized in Algorithm~\ref{alg:solveparllel}.

\begin{algorithm}[t]
\caption{$\textbf{Solver\_Parallel}\left(Ambig,P_m\right)$} \label{alg:solveparllel}
\begin{flushleft}
\textbf{Input:} $Ambig$, $P_m$\\
\textbf{Output}: $\text{STATUS}$, $x_{\text{Sol}}$
\end{flushleft}
\begin{algorithmic}[1]
\STATE The tool runs off-the-shelf solvers such as Z3 or Yices on small-volume ambiguous regions in $Ambig$ in parallel:
    \STATE $\text{STATUS}, x_{\text{Sol}}:=\textbf{Z3/Yices\_Parall}\left(Ambig,P_m\right)$
\RETURN $\text{STATUS}, x_{\text{Sol}}$
\end{algorithmic}
\end{algorithm}


\section{Extension to SMT solving}
We extend the PolyAR solver described in the previous sections to account for combinations of Boolean and Polynomial inequality constraints of the form:
\begin{align}
    \exists& (b_1,\ldots,b_o, x_1,\ldots,x_n) \in \mathbb{B}^o \times \mathbb{R}^n, \nonumber \\
    &\text{subject to:} \nonumber \\
    & p_i(x_1,\ldots,x_n) \le 0, && i = 1,\ldots,m \\
    &x_k \in [\underline{d}_k, \overline{d}_k], && k = 1,\ldots,n\\
    &\varphi_j(b_1,\ldots,b_o) \; \longleftrightarrow \;\texttt{TRUE}, && j = 1,\ldots,r\\
    &b_l \longleftrightarrow \big( p_{l+m}(x_1,\ldots,x_n) \; \le \; 0 \big), && l = 1,\ldots,h \label{eq:hybrid}
\end{align}
where $\varphi_j(b_1,\ldots,b_o)$ is any combinations of Boolean and pseudo-Boolean predicates.

We can create a Satisfiability Modulo Theory (SMT) solver by combining a SAT solver for Boolean and pseudo-Boolean constraints and a theory solver (PolyAR) for interval and polynomial constraints on real numbers by following the lazy SMT paradigm \cite{lazysmt}. The SAT solver solves the combination of Boolean and pseudo-Boolean constraints using the David-Putnam-Logemann-Loveland (DPLL) algorithm and suggests satisfying assignments for the Boolean variables $b$ and thus suggesting which polynomial constraints should jointly satisfied (or unsatisfied). The theory solver (PolyAR) checks the validity of the given assignments and provides an explanation of the conflict, i.e., an \textit{UNSAT certificate}, whenever a conflict is found. Each certificate is a new Boolean constraint that will be used by the SAT solver to prune the search space. 

While in the lazy SMT paradigm, the PolyAR solver needs to be executed multiple times with a different set of polynomial constraints, we modify the PolyAR solver to perform all the abstraction refinement for all the polynomials as a pre-processing step. This eliminates the need to re-compute the same abstraction refinement every time the PolyAR solver is executed. 

Whenever the SAT solver assigns one of the Boolean variables $b_l$ in~\eqref{eq:hybrid} to zero, then the PolyAR solver needs to guarantee that the corresponding polynomial $p_{l+m}$ satisfy $p_{l+m}(x) > 0$ or equivalently $-p_{l+m}(x) \le 0$. To eliminate the need to apply the convex abstraction refinement process for both $p_{l+m}(x)$ and $-p_{l+m}(x)$, the PolyAR solver computes the negative and positive boxes ($\mathcal{B}^N$ and $\mathcal{B}^P$) only for $p_{l+m}(x)$ and flips their usage for $-p_{l+m}(x)$.

\vspace{-3mm}
\section{NUMERICAL RESULTS}
\vspace{-3mm}
In this section, we compare the performance of PolyAR to the state-of-the-art solvers Z3 8.9 and Yices 2.6. The objective of this comparison is to study the performance on:
\begin{itemize}
    \item Problems that appear naturally in parametric controller synthesis. In particular, we focus on the problem of designing stabilizing SOF controllers for LTI systems~\cite{SOFdesign}. 
    \item Problems that appear in non-parametric controller synthesis for non-linear systems. In particular, we focus on the problem of designing a controller for the nonlinear Duffing oscillator~\cite{MPCdesign}.
    \item Additionally, we demonstrate the performance of PolyAR on designing a hybrid switching system; a problem which state-of-the-art tools are incapable of handling.
\end{itemize}
All the experiments were executed on an Intel Core i7 2.6-GHz processor with 16 GB of memory. 
\vspace{-2mm}
\subsection{Static Output Feedback Controller Synthesis for Linear Time Invariant Systems}
\vspace{-2mm}
In this subsection, we assess the scalability of the PolyAR solver compared to state-of-the-art solvers on control synthesis problems. In particular, we consider the problem of synthesizing a parametric controller for the following continuous LTI system:
$$ \dot{x} = A x + B u, \qquad y = C x,$$
where $x \in \mathbb{R}^{n_A}$ is the system state, $u \in \mathbb{R}^{n_B}$ is the system control input, $y \in \mathbb{R}^{n_C}$ is the system output, and the matrices $A \in \mathbb{R}^{n_{A}\times n_{A}}$, $B \in \mathbb{R}^{n_A \times n_B}$ and $C \in \mathbb{R}^{n_{C} \times n_A}$ are the system matrices. We are interested in designing a static output feedback controller of the form:
$$ u =  K y,$$
such that the resulting closed loop system:
$$ \dot{x} = (A + BKC) x,$$
is stable, i.e., the matrix $A+BKC$ is Hurwitz.

We follow the steps detailed in~\cite{SOFdesign} to pose the problem of designing the static output feedback controller as a set of polynomial constraints using the Routh-Hurwitz stability criteria. The Routh-Hurwitz stability criteria result in a set of $n_A$ polynomials in the elements of the controller matrix $K$. We consider five instances of the controller synthesis problem with the following parameters:
\begin{itemize}
    \item \textbf{Example 1}: $n_A = 3, n_B = 4, n_C = 4$ which results in 3 polynomial constraints with 16 variables and max polynomial order of $4$. We restrict the elements of the controller matrix to be inside $[-4,7]$. 
    
    \item \textbf{Example 2}: $n_A = 3, n_B = 5, n_C = 5$ which results in 3 polynomial constraints with 25 variables and max polynomial order of $3$. We restrict the elements of the controller matrix to be inside $[-0.5,1]$.
    
    \item \textbf{Example 3}: $n_A = 2, n_B = 6, n_C = 6$ which results in 2 polynomial constraints with 36 variables and max polynomial order of $3$. We restrict the elements of the controller matrix to be inside $[0,5]$.
    
    \item \textbf{Example 4}: $n_A = 2, n_B = 7, n_C = 7$ which results in 2 polynomial constraints with 49 variables and max polynomial order of $2$. We restrict the elements of the controller matrix to be inside $[-10,0]$.
    
    \item \textbf{Example 5}: $n_A = 5, n_B = 4, n_C = 4$ which results in 5 polynomial constraints with 16 variables and max polynomial order of $4$. We restrict the elements of the controller matrix to be inside $[-4,7]$. In addition, we want to enforce the following controller structure:
    $$k_{21}\times k_{22} \times k_{23} <~0, \quad k_{21}+ k_{22} + k_{23} <-1$$,
    which can be encoded using the additional SMT constraints:
    \begin{align*}
        &b_1 \wedge b_2 \longleftrightarrow \text{True},\\
        & b_1 \rightarrow k_{21}\times k_{22} \times k_{23} < 0,\\
        & b_2 \rightarrow k_{21}+ k_{22} + k_{23} < -1,
    \end{align*}
    where $k_{ij}$ are the elements of the controller matrix $K$.
\end{itemize}
For each of these examples, we generate random system matrices from a zero-mean normal distribution and feed them to four versions of our solver PolyAR:
\begin{itemize}
    \item PolyAR + Z3 (1 thread): This version uses one instance of Z3 to analyze all the ambiguous regions.
    \item PolyAR + Z3 (max threads): This version uses a separate instance of Z3 to analyze each of the ambiguous regions. All Z3 instances are running in parallel.
    \item PolyAR + Yices (1 thread): This version uses one instance of Yices to analyze all the ambiguous regions.
    \item PolyAR + Yices (max thread): This version uses a separate instance of Yices to analyze each of the ambiguous regions. All Yices instances are running in parallel.
\end{itemize}
We compare the execution times of these four solvers with Z3 8.9 and Yices 2.6. Table~1 shows the execution time for all the solvers. As evident by the results in Table~1, off-the-shelf solvers are incapable of solving all the five examples and they time out after one hour. On the other hand, and thanks to the abstraction refinement process, the PolyAR solver is able to solve all the instances in a few seconds, leading to $240X$ speed up in the total execution time in the PolyAR+Yices (max threads) case, evidence of the scalability of the proposed approach.

\begin{table}[t!]
Table 1: Experiment results for SOF design. The timeout is set by $3600~s$.

\begin{adjustbox}{width=\columnwidth,center}
\begin{tabular}{|c|c|c|c|c|c|c|}
    \hline
     \multirow{3}{*}{Example} & \multicolumn{6}{c|}{Times (seconds)} \\
    \cline{2-7}
    & Z3 8.9  & Yices 2.6 & PolyAR+Z3 & PolyAR+Z3 & PolyAR+Yices & PolyAR+Yices\\
    & & & (1 thread) & (max threads) & (1 thread) & (max threads)\\
    \hline
    \hline
    1  & \textit{timeout} & \textit{timeout} & \textit{timeout} & $7.552$ & $\mathbf{2.405}$ & $2.442$\\
   \hline
   2   & \textit{timeout} & \textit{timeout} & $83.776$ & $114.453$ & $timeout$ & $\mathbf{3.766}$ \\
    \hline 
   3  & \textit{timeout} & \textit{timeout} & $23.551$ & $23.970$ & \textit{timeout} & $\mathbf{8.725}$ \\
     \hline 
   4  & \textit{timeout} & \textit{timeout} & $0.718$ & $0.729$ & $\mathbf{0.416}$ & $0.432$\\
    \hline 
   5  & \textit{timeout} & \textit{timeout} & $3.636$ & $3.768$ & $0.621$ & $\mathbf{0.498}$  \\
      \hline 
     \hline
    \# Problems   & $0$ & $0$ & $4$ & $5$ & $3$ & $\mathbf{5}$ \\
     Solved & & & & & & \\
\hline
    Total Time  & \textit{timeout} & \textit{timeout} & $111.681$ & $150.472$ & $3.442$ & $\mathbf{15.863}$   \\
    (seconds)  & & & & & &\\
    \hline
\end{tabular}
\end{adjustbox}
\vspace{-2mm}
\end{table}

In the following, we give the stabilizing controller matrices $K_{1}$, $K_{2}$, $K_{3}$, $K_{4}$, $K_{5}$, and the two Boolean variables $b_1$ and $b_2$ that PolyAR (Yices) returned for Examples 1, 2, 3, 4, and 5:     
\begin{align*}
K_1&=\begin{bmatrix}
-4 & -2 & -2 & 1 \\
1 & 1 & 1 & 1 \\
1 & 1 & 1 & 1 \\
1 & 4 & 1 & 1  
\end{bmatrix}, K_2=\begin{bmatrix}
0 & 1 & 0 & 0 & 0\\
0 & 0 & 0 & 0 & 0\\
0 & 0 & 0 & 1 & 0\\
0 & 0 & 0 & 0 & 0\\
0 & 0 & 0 & 0 & 0
\end{bmatrix},\\
K_3&=\begin{bmatrix}
1 & 1 & 3 & 3 & 2 & 2 \\
 2 & 2 & 2 & 2 & 2 & 2 \\
 2 & 0 & 0 & 2 & 2 & 5 \\
 0 & 0 & 2 & 2 & 0 & 2 \\
 2 & 2 & 0 & 0 & 2 & 5 \\
 5 & 2 & 2 & 2 & 2 & 2  
\end{bmatrix}, K_4=\begin{cases}
    -8, &  i=1,j=1,\\  -5,         &   2\leq i,j\leq 7,
\end{cases}\\
K_5&=\begin{bmatrix}
-2 & 4 & -2 & 7 \\
1 & -3.5 & 1 & 1 \\
5 & 6 & 1 & 1 \\
1 & 6.99 & 1 & 1  
\end{bmatrix},~b_1=\text{True},~b_2=\text{True}.
\end{align*}
It is easily to note that the solutions given by PolyAR (Yices) satisfies the two Boolean constraints, i.e., $k_{21}\times k_{22} \times k_{23}=-3.5 <~0$ and $k_{21}+ k_{22} + k_{23}=-1.5 <-1$. 

In conclusion, PolyAR+Yices (max thread) solver outperforms all the other solvers due to the effectiveness of Yices in reasoning about problems with small volumes.

\subsection{Non-Linear Controller Design for a Duffing Oscillator}
In this subsection, we assess the scalability of PolyAR solver compared to state-of-the-art solvers on synthesizing a non-parametric controller for a Duffing oscillator reported by~\cite{MPCdesign}. The dynamics of the oscillator is given by the higher-order differential equation:
\begin{align}\label{duff1eq}
   y^{(n)}\!\!\left(t\right)\!+\!\cdots\!+\!y^{(2)}\!\!\left(t\right)\!+\!2\zeta y^{(1)}\!\!\left(t\right) \!+\! y\!\left(t\right) \!+\! y\!\left(t\right)^3 \!\!=\!\! u\left(t\right),  
\end{align}
where $y \in \mathbb{R}$ is the continuous state variable and $u \in \mathbb{R}$ is the control input. The parameter $\zeta$ is the damping coefficient. The objective of the control is to regulate the state to the origin. To derive the discrete-time model, forward difference approximation is used (with sampling period of $h=0.05$ time units). The resulting state space model with discrete state vector $x=[x_1,x_2,\cdots,x_n]^T=[y,y^{(1)},\cdots,y^{(n-1)}]^T \in \mathbb{R}^{n-1}$ and input $u \in \mathbb{R}$ is:  
\begin{align}\label{duff2eq}
\begin{bmatrix}
x_1\\
x_2\\
\vdots\\
x_n
\end{bmatrix}^{+}&\!\!\!\!\!\!=\!\!\begin{bmatrix}
\!1\! & \!\!\!h\! & \!\!\!0\! & \!\!\!\cdots\! & \!\!\!\cdots\! & \!\!\!0\!\\
\!0\! & \!\!\!1\! & \!\!\!h\! & 0 & \!\!\!\cdots\! & \!\!\!0\!\\
\!\vdots\! & \!\!\!\vdots\! & \!\!\!\vdots\! & \!\!\!\vdots\! & \!\!\!\vdots\! & \!\!\!\vdots\!\\
\!-h\! & \!\!\!-2\zeta h\! & \!\!\!-h\! & \!\!\!\hdots\! & \!\!\!-h\! & \!\!\!1\!\!-\!\!h\!
\end{bmatrix}\!\!\!\!
\begin{bmatrix}
x_1 \\
x_2 \\
\vdots\\
x_2
\end{bmatrix}
\!\!\!+\!\!\!\begin{bmatrix}
0 \\
0 \\
\vdots\\
h
\end{bmatrix}\!\!u\!+\!\!\!\begin{bmatrix}
0 \\
0 \\
\vdots\\
\!\!-hx_1^3\!
\end{bmatrix}.
\end{align}
The previous equation 
is written in the form of $x\left(k+1\right)=Ax\left(k\right)+Bu\left(k\right)+E\left(x\right)$, which includes a nonlinear term $E\left(x\right)=\begin{bmatrix} 0, \cdots, -hx_1^3\left(k\right) \end{bmatrix}^T$. Our objective is to design a non-parametric controller. To that end, we encode the controller as the solution of a feasibility problem of several constraints that capture the system dynamics, state/input constraints, and stability constraints as discussed below.

First, to enforce the stability of the resulting non-parametric controller, we consider the candidate quadratic Lyapunov function $V\left(x\right)=x^TPx$ with the symmetric positive definite matrix $P$ is a solution of the discrete-time Lyapunov equation $APA^T+P+Q=0$ and is a positive definite matrix. Thanks to the fact that $E(x)$ satisfies $\text{lim}_{\norm{x} \rightarrow 0}\frac{\norm{E\left(x\right)}}{\norm{x}}=0$ along with the Lyapunov's indirect method in~\cite{khalil}, one can directly conclude that $V(x)$ is indeed a Lyapunov function.
For simplicity, we pick $Q=I_n$, where $I_n$ is the identity matrix of size $n$. 

Moreover, to ensure the smoothness of the resulting controller signals, we add additional filters in the form of high order polynomial $L(x,u) \leq 0$. 
In addition, we consider the state-constraints of the form $\norm{x\left(k\right)}_{\infty} \leq 0.6$.

The final non-parametric controller is then encoded as the solution of the following feasibility problem:
\begin{align}\label{mpc4}
&\exists x_1(k),\ldots x_n(k), x_1(k+1),\ldots x_n(k+1), u(k)  \nonumber\\
&\text{subject to}: \nonumber\\
&\qquad x\left(k+1\right)=Ax\left(k\right)+Bu\left(k\right)+E\left(x\right), \nonumber\\
&\qquad V\left(x\left(k+1\right)\right)-V\left(x\left(k\right)\right) \le -\epsilon, \nonumber\\
&\qquad L\left(x\left(k\right), u(k)\right) \leq 0, \nonumber \\
&\qquad \norm{x\left(k\right)}_{\infty} \leq 0.6.
\end{align}
Since the PolyAR solver only handles polynomial inequalities, hence, we transform the equality constraint 
$x\left(k+1\right)=Ax\left(k\right)+Bu\left(k\right)+E\left(x\right) $ above into two inequalities $x\left(k+1\right)-Ax\left(k\right)+Bu\left(k\right)+E\left(x\right) \leq\epsilon~\wedge~x\left(k+1\right)-Ax\left(k\right)+Bu\left(k\right)+E\left(x\right) \geq -\epsilon$, where $\epsilon \in \mathbb{R}$ is a small value. 

We consider three instances of the controller synthesis problem for the Duffing oscillator with the following parameters:
\begin{itemize}
    \item $n=2$, $\zeta=0.3$, $x\left(0\right)=[0.4,0.1]^T$, $L\left(x\left(k\right),u\left(k\right)\right)=x_1^{11}\left(k\right)+x_2^{11}\left(k\right)-u^{10}\left(k\right)$, which results in $6$ polynomial constraints with $3$ variables and max polynomial order of $11$.
    \item $n=3$, $\zeta=1.0$, $x\left(0\right)=[0.1,0.1,0.1]^T$, $L\left(x\left(k\right),u\left(k\right)\right)=x_1^{5}\left(k\right)+x_2^{5}\left(k\right)+x_3^{5}\left(k\right)+u^{5}\left(k\right)$, which results in $8$ polynomial constraints with $4$ variables and max polynomial order of $5$.
    \item $n=4$, $\zeta=1.75$, $x\left(0\right)=[0.1,0.1,0.01, 0.1]^T$, $L\left(x\left(k\right),u\left(k\right)\right)=x_1^{4}\left(k\right)+x_2^{4}\left(k\right)+x_3^{4}\left(k\right)+x_4^{4}\left(k\right)-u^{4}\left(k\right)$, which results in $10$ polynomial constraints with $5$ variables and max polynomial order of $4$.
\end{itemize}

\begin{figure}
\resizebox{.45\textwidth}{!}{
    \centering
    \begin{tabular}{ c | c | c |}
         $n$ & State Space & Execution Time Evolution over time\\\hline
         2 & 
         \raisebox{-0.5\totalheight}{\includegraphics[width=0.6\columnwidth,trim=8mm 0 15mm 0, clip]{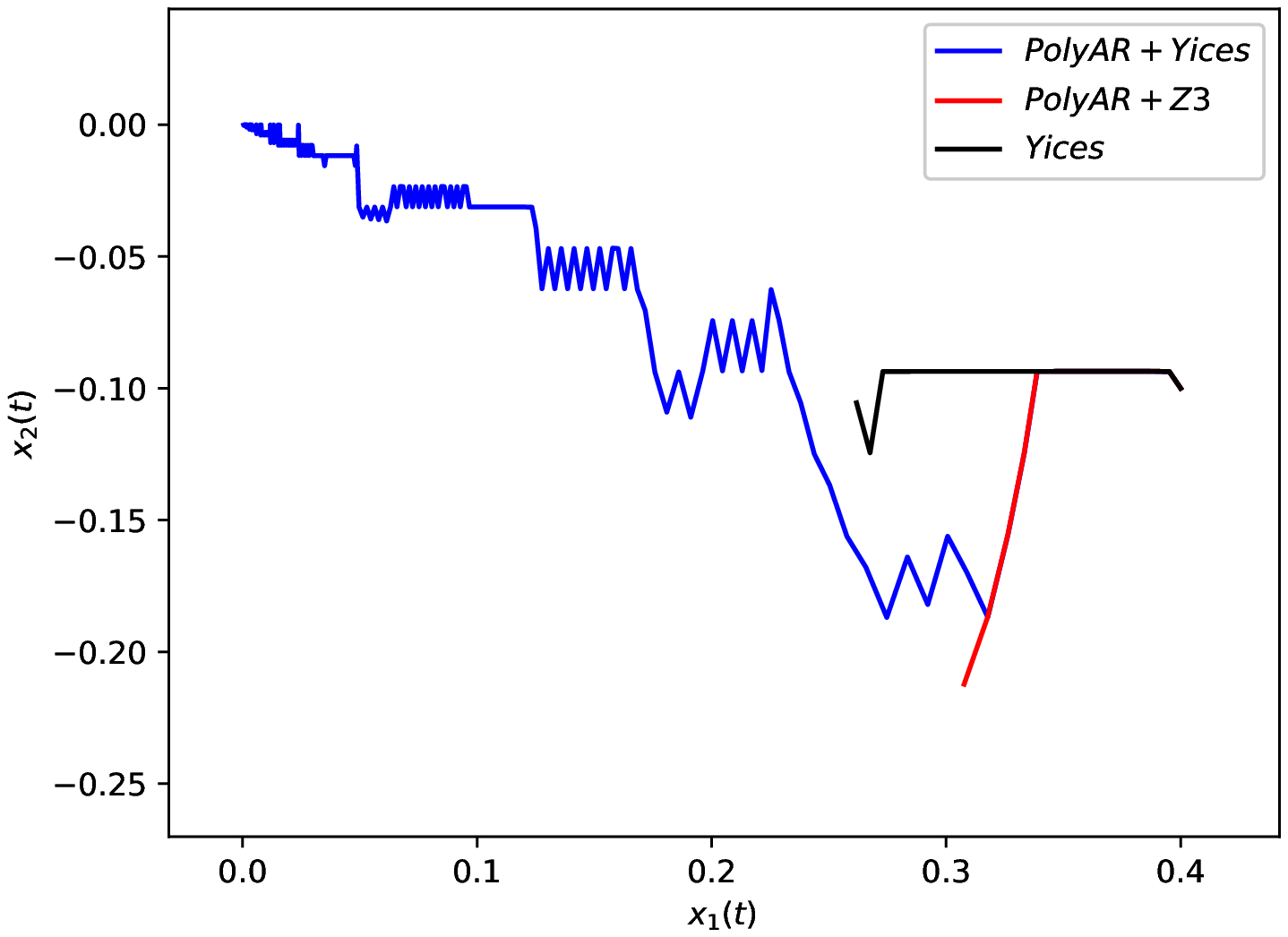}} &
         \raisebox{-0.5\totalheight}{\includegraphics[width=0.6\columnwidth,trim=8mm 0 15mm 0, clip]{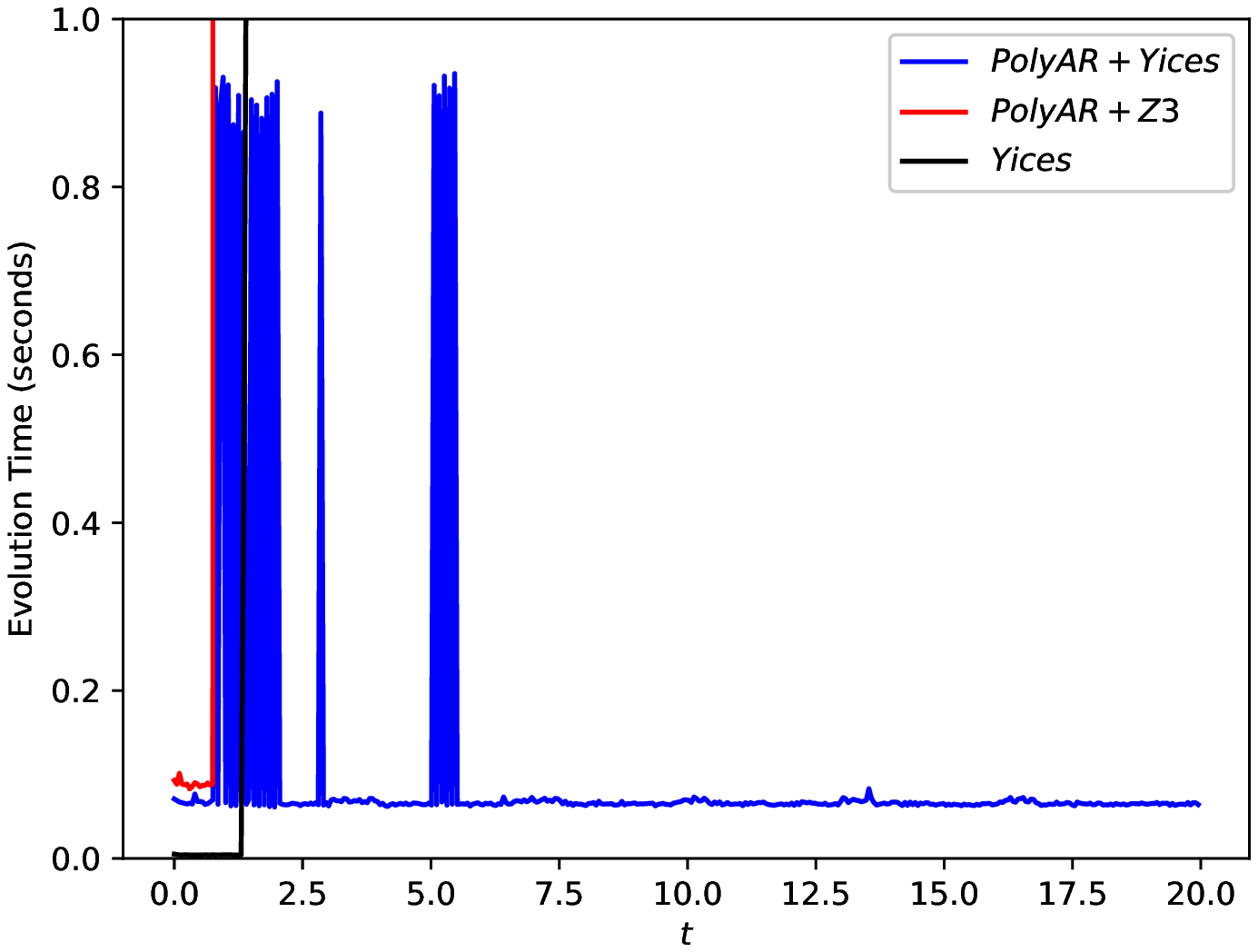}} \\ \hline
          3 & 
         \raisebox{-0.5\totalheight}{\includegraphics[width=0.6\columnwidth,trim=8mm 0 15mm 0, clip]{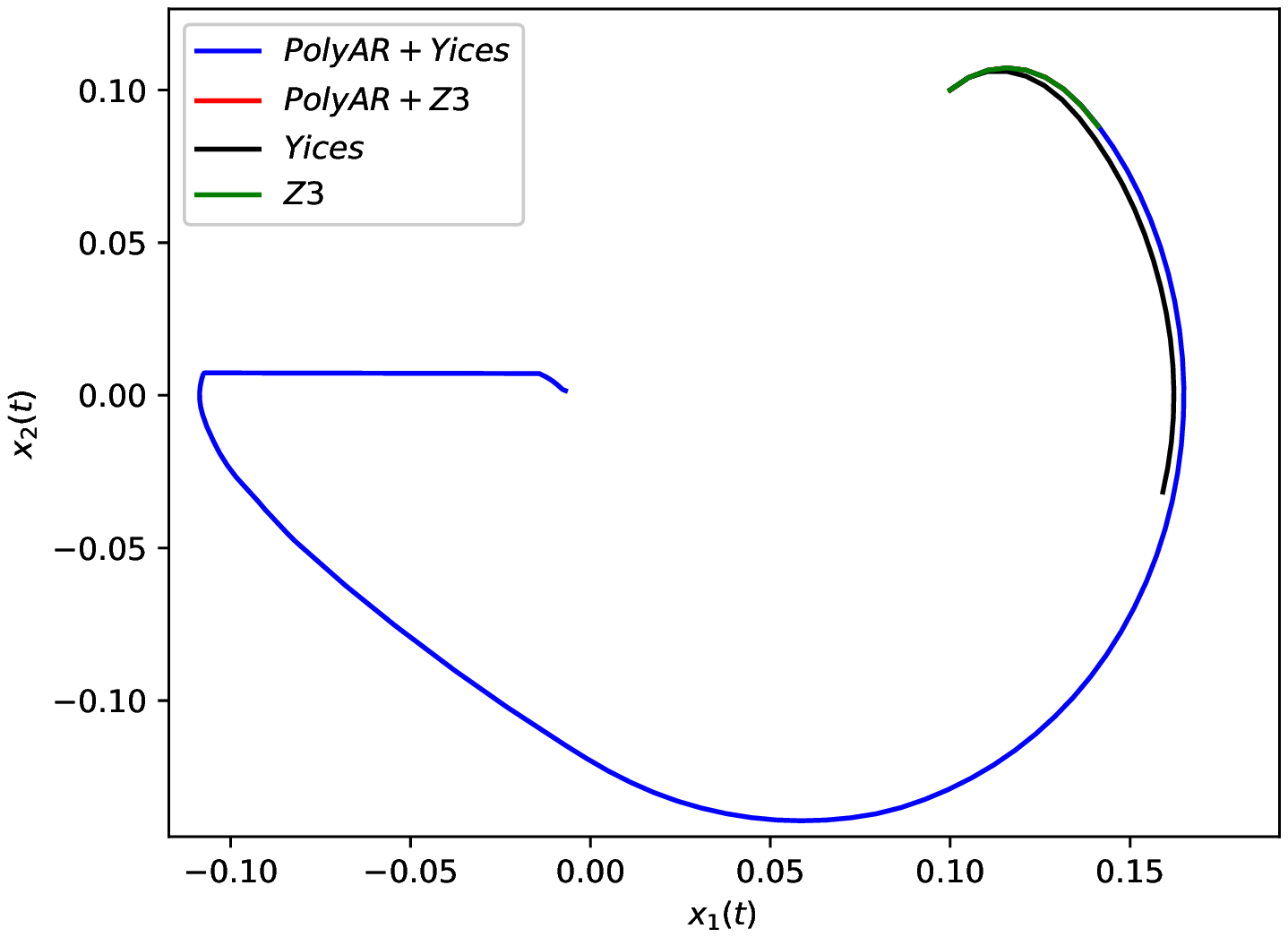}} &
         \raisebox{-0.5\totalheight}{\includegraphics[width=0.6\columnwidth,trim=8mm 0 15mm 0, clip]{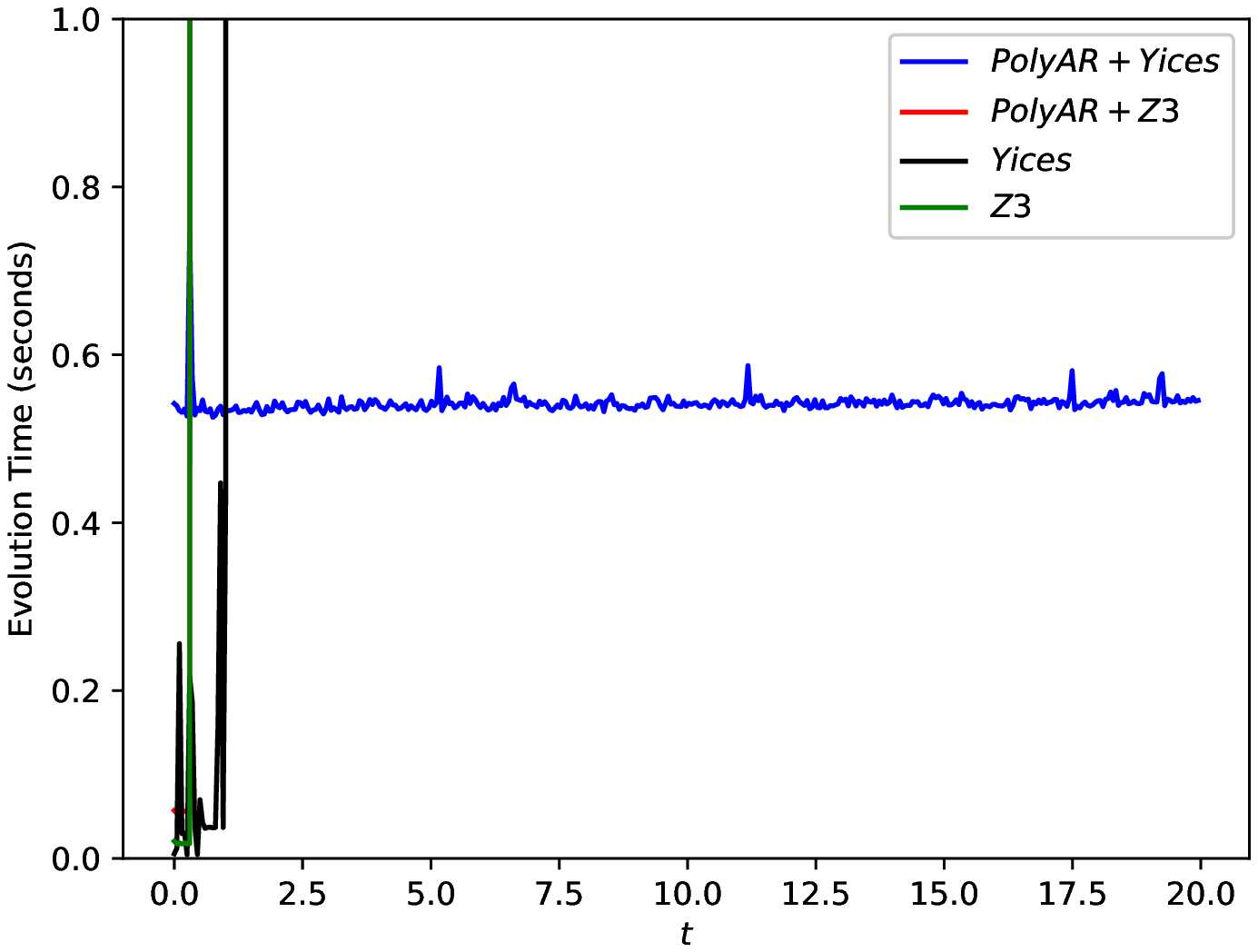}} \\ \hline
          4 & 
         \raisebox{-0.5\totalheight}{\includegraphics[width=0.6\columnwidth,trim=8mm 0 15mm 0, clip]{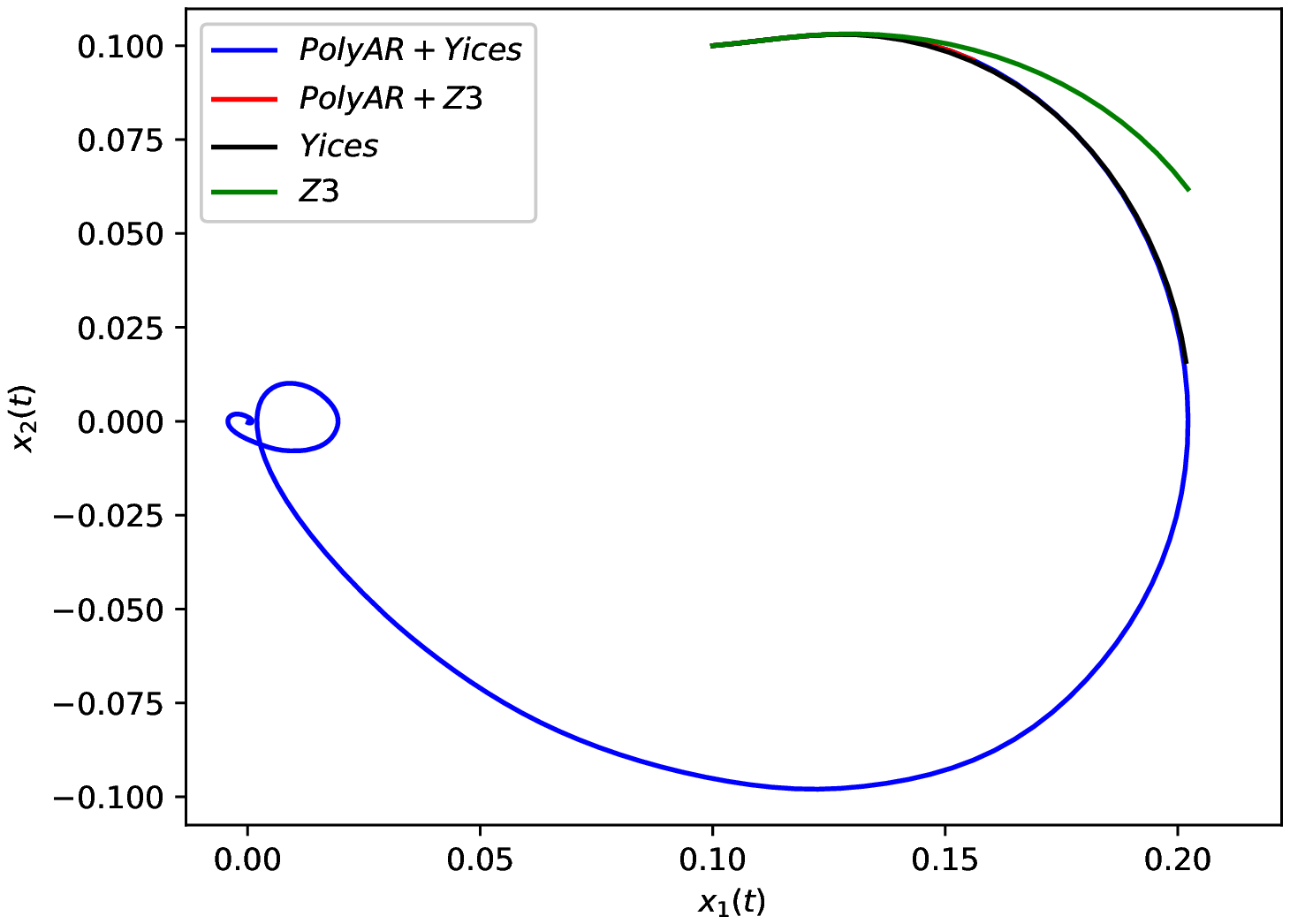}} &
         \raisebox{-0.5\totalheight}{\includegraphics[width=0.6\columnwidth,trim=8mm 0 15mm 0, clip]{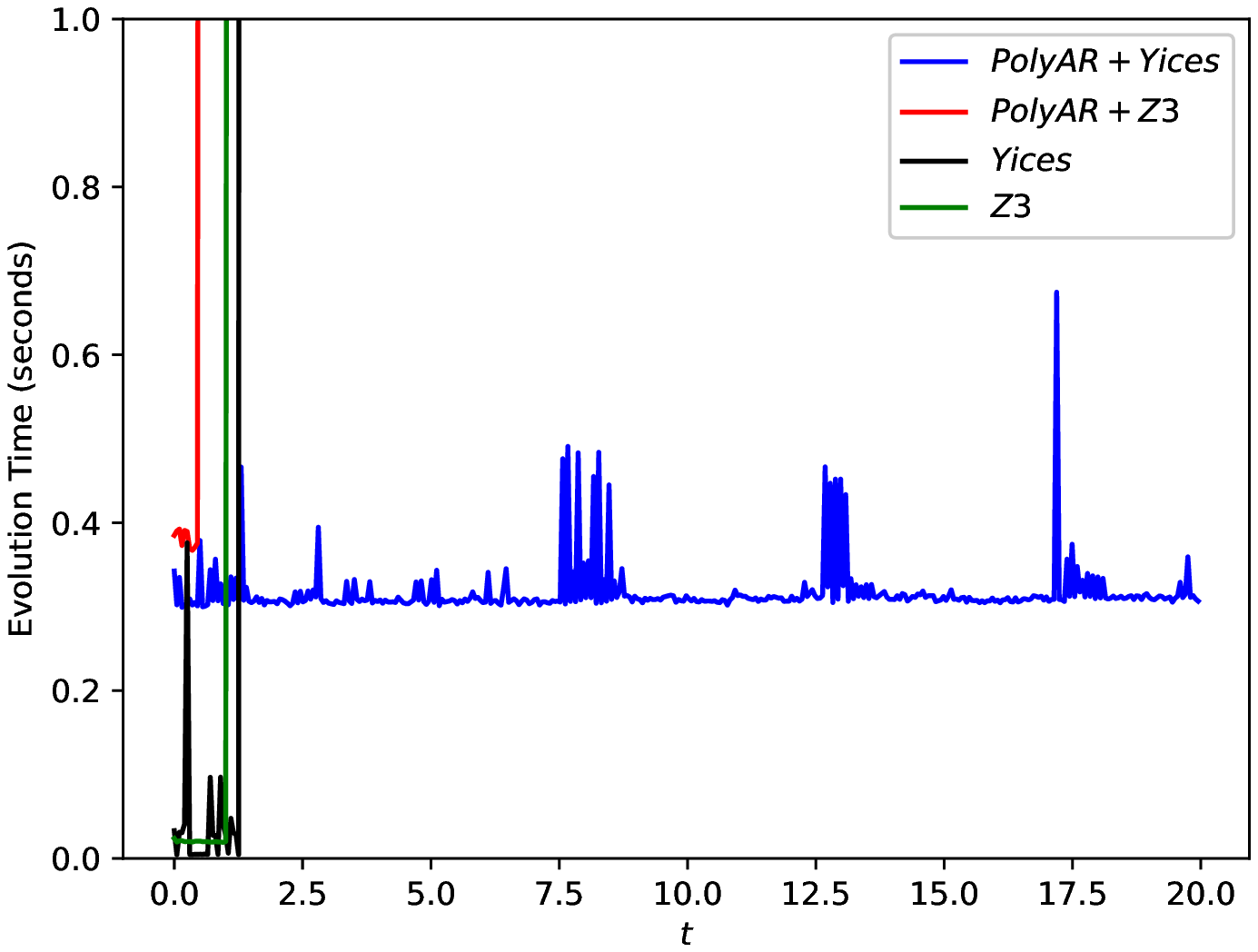}} \\ \hline
    \end{tabular}
    }
    \caption{\small{Results of controlling the Duffing oscillator with different $n$ (left) evolution of the states $x_1(k)$ and $x_2(k)$ for the solvers in the state-space, (right) evolution of the execution time of solvers during the $20$ seconds. The timeout is equal to $1s$. Trajectories are truncated once the solver exceeds the timeout limit.}}
    \label{fig:duff}
    \vspace{-1mm}
\end{figure}

We feed the resultant polynomial inequality constraint to PolyAR$+$Yices, PolyAR$+$Z3, Yices, and Z3. We solve the feasibility problem for $n=2$, $n=3$, and $n=4$. We set the timeout to be $1s$. Figure~\ref{fig:duff} (left) shows the state-space evolution of the controlled Duffing oscillator for different solvers for number of variables $n$ of $2,3,$ and $4$. Figure~\ref{fig:duff} (right) shows the evolution of the execution time of the solvers during the $20$ seconds. As it can be seen from Fig. 4, our solver PolyAR$+$ Yices succeeded to find a control input $u$ that regulates the state to the origin for all $n$. However, off-the-shelf solvers are incapable of solving all the three instances and they early time out after one second out of the simulated $20~seconds$.

\begin{figure}[!ht]
		\centering
		\includegraphics[width=0.85\columnwidth]{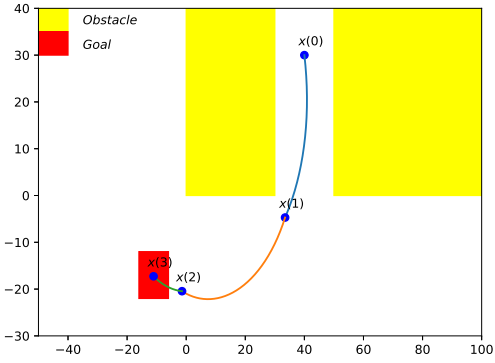}
		\caption{\small{The trajectory  that starts from an initial state $x\left(0\right)=[40,30]^T$ and reaching a final state $x\left(3\right) \in Goal$ while avoiding the obstacles. The goal and the obstacles are represented with a red and yellow rectangle, respectively.}}
		\label{expmatr}
\end{figure}

\subsection{Designing Switching Signals for Continuous-Time Linear Switching Systems}
\vspace{-3mm}
In this subsection, we show how to use the PolyAR solver to successfully design a controller for a continuous-time linear switching system. In particular, we consider the following switching dynamics:
$$ \dot{x} = A_{\sigma(t)} x, \qquad \sigma(t) = \{1,2,3\},$$
with $x(t) \in \mathcal{X} \subset \mathbb{R}^2$ is the system state at time $t$ and the matrices $A_1$, $A_2$, and $A_3 \in \mathbb{R}^{2\times2}$ represents three modes for the switching system. Consider the state space in Figure~\ref{expmatr}. The objective is to design a switching signal $\sigma(t)$ that can steer the state of the system to the goal set $Goal~\subset~\mathcal{X}$ while avoiding entering the obstacle set $Obstacle~\subset\mathcal{X}$. For simplicity, we confine our attention to step-wise switching signals $\sigma(t)$. That is, we assume the switching signal $\sigma(t)$ will be constant for some amount of time $t_1, t_2, \ldots, t_L$. Our objective is then to design the switching times and the associated system mode that leads to the satisfaction of the reach-avoid specifications. To that end, we define a set of Boolean variables $b_{ij}$ such that $b_{ij}$ is equal to 1 whenever the $j$th mode is active during $t_i$. Given the initial condition of the system $x(0)$, we can use these Boolean variables to encode the problem of designing the switching signal as the following SMT constraints:
\begin{align}\label{switch1}
&\exists b_{11},\ldots,b_{13},\ldots,b_{L1}, \ldots, b_{L3}, x(1), \ldots, x(L), t_1, \ldots, t_L \nonumber\\
&\text{subject to:} \nonumber\\
&\qquad b_{11}\rightarrow x\left(1\right)=\exp{\left(A_1t_1\right)}~x\left(0\right),\nonumber\\
&\qquad b_{12}\rightarrow x\left(1\right)=\exp{\left(A_2t_1\right)}~x\left(0\right),\nonumber\\
&\qquad b_{13}\rightarrow x\left(1\right)=\exp{\left(A_3t_1\right)}~x\left(0\right),\nonumber\\
&\qquad b_{11}+b_{12}+b_{13}=1,\nonumber\\
&\qquad \qquad\qquad \vdots\nonumber\\
&\qquad b_{L1}\rightarrow x\left(L\right)=\exp{\left(A_1t_L\right)}~x\left(L-1\right),\nonumber\\
&\qquad b_{L2}\rightarrow x\left(L\right)=\exp{\left(A_2t_L\right)}~x\left(L-1\right),\nonumber\\
&\qquad b_{L3}\rightarrow x\left(L\right)=\exp{\left(A_3t_L\right)}~x\left(L-1\right),\nonumber\\
&\qquad b_{L1}+b_{L2}+b_{L3}=1,\nonumber\\
&\qquad x\left(1\right),\ldots,x\left(L-1\right) \notin Obstacle,\nonumber\\
&\qquad x\left(L\right) \in Goal,
\end{align}
where the constraint $b_{i1} + b_{i2} + b_{i3} = 1$ is a pseudo-Boolean constraint that enforces the consistency between the Boolean variables such that only one of the three modes $A_1,A_2$ and $A_3$ can be selected during the period $t_{j-1} < t \le t_j$.
Since PolyAR solver only handles polynomial inequalities, we approximate the exponential matrix $\exp\left(A_it_j\right)\approx I_2+t_jA_i+\frac{t_j^2A_i^2}{2}+\frac{t_j^3A_i^3}{6}$, $i,j=1,\cdots,L$, where $I_2$ is the identity matrix of size $2$ and $A_i^n=A_i\times \cdots \times A_i$. Furthermore, we transform the equality $x\left(i\right)=\exp{\left(A_it_j\right)}~x\left(i-1\right)$, $i,j=1,\cdots,L$, into two inequalities $x\left(i\right)-\exp{\left(A_it_j\right)}~x\left(i-1\right) \leq \epsilon \wedge x\left(i\right)-\exp{\left(A_it_j\right)}~x\left(i-1\right) \geq -\epsilon$, where $\epsilon \in \mathbb{R}$ is a small value.

In our experiments, we pick the horizon $L = 3$ and the modes $A_1=\begin{bmatrix}-1, 2\\-2,-2 \end{bmatrix}$, $A_2=\begin{bmatrix}-1, 3\\-3,-1 \end{bmatrix}$, and $A_3=\begin{bmatrix}0, 2\\-2,0 \end{bmatrix}$, and we start with an initial state $x\left(0\right)=[40,30]^T$. We restrict the states to be inside $\mathcal{X} \in [-100,100]$. We feed the resultant polynomial inequality constraint to PolyAR+Yices, and as it can be seen from Fig.\ref{expmatr}, the solver succeeded to find the right modes ($b_{11}=b_{22}=b_{33}=1$) and the necessary times $t_1=0.391s$, $t_2=0.5s$, and $t_3=0.25s$ that ensures that $x\left(3\right)$ reaches a $Goal$ while the intermediate states $x\left(2\right),x\left(1\right)$ avoid $Obstacles$. In addition, we remark that the trajectory between $x\left(0\right)$ and $x\left(2\right)$ is making its way to the equilibrium point $[0,0]^T$. This is explained by the fact that the matrices $A_1$ and $A_2$ are stables. Our solver computes the necessary time $t_3$ that ensures that the final state $x\left(3\right) \in~Goal$ and does not converge to the equilibrium point.
\vspace{-3mm}
\bibliography{main}             

\begin{thebibliography}{13}
\providecommand{\natexlab}[1]{#1}
\providecommand{\url}[1]{\texttt{#1}}
\providecommand{\urlprefix}{URL }
\expandafter\ifx\csname urlstyle\endcsname\relax
  \providecommand{\doi}[1]{doi:\discretionary{}{}{}#1}\else
  \providecommand{\doi}{doi:\discretionary{}{}{}\begingroup
  \urlstyle{rm}\Url}\fi

\bibitem[{{Bahavarnia} et~al.(2020){Bahavarnia}, {Shoukry}, and
  {Martins}}]{SOFdesign}
{Bahavarnia}, M., {Shoukry}, Y., and {Martins}, N.C. (2020).
\newblock {Controller Synthesis subject to Logical and Structural Constraints:
  A Satisfiability Modulo Theories (SMT) Approach}.
\newblock In \emph{2020 American Control Conference (ACC)}, 5281--5286.
\newblock \doi{10.23919/ACC45564.2020.9147460}.

\bibitem[{Barrett and Tinelli(2018)}]{lazysmt}
Barrett, C. and Tinelli, C. (2018).
\newblock {Satisfiability modulo theories}.
\newblock In \emph{Handbook of Model Checking}, 305--343. Springer.

\bibitem[{{Bauer} et~al.(2007){Bauer}, {Pister}, and {Tautschnig}}]{ABsolver}
{Bauer}, A., {Pister}, M., and {Tautschnig}, M. (2007).
\newblock {Tool-support for the analysis of hybrid systems and models}.
\newblock In \emph{2007 Design, Automation Test in Europe Conference
  Exhibition}, 1--6.

\bibitem[{Behroozi(2019)}]{InscribRec}
Behroozi, M. (2019).
\newblock {Largest Inscribed Rectangles in Geometric Convex Sets}.
\newblock \emph{CoRR}, abs/1905.13246.

\bibitem[{Collins(1975)}]{collins}
Collins, G.E. (1975).
\newblock {Quantifier elimination for real closed fields by cylindrical
  algebraic decomposition}.
\newblock In \emph{Automata Theory and Formal Languages 2nd GI Conference
  Kaisers lautern, May 20--23, 1975}, 134--183. Springer.

\bibitem[{De~Moura and Bj{\o}rner(2008)}]{Z3}
De~Moura, L. and Bj{\o}rner, N. (2008).
\newblock {Z3: An efficient SMT solver}.
\newblock In \emph{International conference on Tools and Algorithms for the
  Construction and Analysis of Systems}, 337--340.

\bibitem[{England and Davenport(2016)}]{complexityproblem1}
England, M. and Davenport, J.H. (2016).
\newblock {The complexity of cylindrical algebraic decomposition with respect
  to polynomial degree}.
\newblock In \emph{International Workshop on Computer Algebra in Scientific
  Computing}, 172--192. Springer.

\bibitem[{Ferreau et~al.(2016)Ferreau, Alm{\'e}r, Peyrl, Jerez, and
  Domahidi}]{surveuNumTool1}
Ferreau, H.J., Alm{\'e}r, S., Peyrl, H., Jerez, J.L., and Domahidi, A. (2016).
\newblock {Survey of industrial applications of embedded model predictive
  control}.
\newblock In \emph{2016 European Control Conference (ECC)}, 601--601. IEEE.

\bibitem[{Fotiou et~al.(2006)Fotiou, Rostalski, Parrilo, and
  Morari}]{MPCdesign}
Fotiou, I.A., Rostalski, P., Parrilo, P.A., and Morari, M. (2006).
\newblock {Parametric optimization and optimal control using algebraic geometry
  methods}.
\newblock \emph{International Journal of Control}, 79(11), 1340--1358.

\bibitem[{Hong(1990)}]{Hong}
Hong, H. (1990).
\newblock {An Improvement of the Projection Operator in Cylindrical Algebraic
  Decomposition}.
\newblock In \emph{Proceedings of the International Symposium on Symbolic and
  Algebraic Computation}, ISSAC '90, 261–264. Association for Computing
  Machinery, New York, NY, USA.
\newblock \doi{10.1145/96877.96943}.
\newblock \urlprefix\url{https://doi.org/10.1145/96877.96943}.

\bibitem[{Khalil(2002)}]{khalil}
Khalil, H.K. (2002).
\newblock \emph{{Nonlinear systems; 3rd ed.}}
\newblock Prentice-Hall, Upper Saddle River, NJ.

\bibitem[{McCallum(1998)}]{McCalum}
McCallum, S. (1998).
\newblock An improved projection operation for cylindrical algebraic
  decomposition.
\newblock In \emph{Quantifier Elimination and Cylindrical Algebraic
  Decomposition}, 242--268. Springer.

\bibitem[{Rao(2009)}]{surveuNumTool2}
Rao, A.V. (2009).
\newblock {A survey of numerical methods for optimal control}.
\newblock \emph{Advances in the Astronautical Sciences}, 135(1), 497--528.

\end{thebibliography}
\end{document}